\documentclass{article}
\usepackage[utf8]{inputenc}

\usepackage{amsmath,amssymb,amsfonts,amsthm}
\usepackage{comment}

\usepackage[left=15mm, top=20mm, right=15mm, bottom=30mm, nohead]{geometry} 

\title{Matrix Models, Integral Polyhedra and Toric Geometry}
\author{Aleksey Andreev\thanks{andreev.av@phystech.edu}}
\date{}

\begin{document}

\maketitle
\vspace{-4.2cm}
\hfill MIPT/TH-20/22

\hfill ITEP/TH-23/22

\hfill IITP/TH-22/22

\vspace{2.7cm}

\setcounter{page}{1}
\begin{center}\vspace{-1cm}
    {\small {\it Moscow Institute of Physics and Technology, Dolgoprudny 141701, Russia }}\\
	{\small {\it NRC, ``Kurchatov Institute'', Moscow 123182, Russia}}\footnote{former Institute for Theoretical and Experimental Physics, Moscow 117218, Russia}
\end{center}

\begin{center}
    ABSTRACT
\end{center}

{\footnotesize
We propose to take a look at a new approach to the study of integral polyhedra. The main idea is to give an integral representation, or matrix model representation, for the key combinatorial characteristics of integral polytopes. Based on the well-known geometric interpretations of matrix model digram techniques, we construct a new model that enumerates triangulations, subdivisions, and numbers of integral points of integral polygons. This approach allows us to look at their combinatorics from a new perspective, motivated by knowledge about matrix models and their integrability. We show how analogs of Virasoro constraints appear in the resulting model. Moreover, we make a natural generalization of this matrix model to the case of polytopes of an arbitrary dimension, considering already a tensor model. We also obtain an analogue of Virasoro constraints for it and discuss their role in the solvability of these models. The deep connection between the geometry of convex polyhedra and toric geometry is the main reference point in the construction of these models. We present considerations on specific ways of applying this approach to the description of Batyrev's mirror pairs. All this allows us to formulate many interesting directions in the study of the connection between matrix/tensor models and the geometry of toric varieties.}

\renewcommand*\contentsname{Contents}

\tableofcontents
\newpage

\section{Introduction}

In this paper, we propose to look at the combinatorics of convex polyhedra and toric geometry from a new perspective. Namely, we propose to apply the method of matrix models to the study of combinatorial problems about integral polyhedra, which in turn are closely related to problems of toric geometry. It seems to us that such an approach can be very useful in the future both from the point of view of toric geometry and from the point of view of matrix models.

Matrix models are an effective method for solving various problems in combinatorics and enumerative geometry. One of the most famous examples is Kontsevich matrix model \cite{Kontsevich_1992}, which enumerates integrals of $\psi$-classes of moduli spaces of complex curves. The distinctive feature of matrix models is the reach variety of analysis methods since they are always finite-dimensional well-defined integrals \cite{Morozov_1994}. It is often occur that matrix models turn out to be solutions of some integrable hierarchy, which makes it possible to ``solve it completely'', which means obtaining exact answers for all (invariant) correlators in this model \cite{mironov2011complete,mironov2017complete,mironov2020cut}. For example, the Kontsevich model satisfies the equations of the KdV hierarchy, and AMM/EO topological recursion \cite{alexandrov2005solving,eynard2007invariants} allows one to write recursive relations on correlators, which are exactly integrals of $\psi$-classes over the moduli space. This example is far from the only one, there are many other problems solved by matrix models (or in which topological recursion can be found) \cite{Morozov_1994,kazarian2015combinatorial,mironov2011complete,kazarian2007algebro,eynard2009topological,eynard2011invariants}. Inspired by many of these examples, we try to apply this approach to describe the combinatorics of convex integral polyhedra.

The combinatorics of convex integral polyhedra is amusing in itself, but we are interested in it because of the well-known connection between the combinatorics of polyhedra and the geometry of toric varieties \cite{fulton1993introduction,buchstaber,audin2012topology}. In particular, we are especially interested in further applying the technique of matrix models to problems of geometry of Batyrev's mirror pairs \cite{batyrev1994dual}, which correspond to reflexive convex integral polyhedra.

In section \ref{section_preliminaries} we give the basic information about integral polyhedra that we need and a brief description of some techniques of matrix models using the Hermitian one as an example.

In this paper, however, we do not discuss applications to Batyrev's theory yet, but we concentrate on the first natural problem in this direction, which is to construct examples of matrix models that describe in some sense the combinatorics of convex integral polyhedra. More specifically, we would like to build matrix models that:
\begin{itemize}
    \item compute the number of integral points in an integral polyhedron,
    \item enumerate triangulations (subdivisions into polytopes) of integral polyhedra
\end{itemize}
The first one, although simple, is quite interesting, because the Hodge numbers of Calabi-Yau hypersurfaces in Batyrev's construction are expressed in terms of the numbers of integral points inside the (sub)polyhedron. The second is directly related to the combinatorics of secondary polytope introduced by Gelfand-Kapranov-Zelevinsky \cite{}.

In order to build a matrix model, we use a simple observation about known models with a cubic potential: it's diagram technique enumerate various triangulations of surfaces with certain weights. The simplest example of such a matrix model is the Hermitian matrix model with a cubic potential whose partition function has the form:

\begin{equation}
    Z=\int \mathcal{D}H \exp\left\{-\frac{1}{2}Tr(H^2)+t_3Tr(H^3)\right\}
\end{equation}

\noindent where the integration is over all Hermitian matrices and the expansion in a series in $t_3$ under the integral is implied. Due to Wick's rule, the terms of the series in $t_3$ enumerate various triangulations of Riemann surfaces of various genera with a fixed number of vertices. There are also more complicated examples of such matrix models, where the correlators can be given the appropriate geometric meaning: the Kontsevich model, the Kazakov model, Hurwitz numbers model and so on.
Based on similar considerations, we construct a matrix model for triangulating two-dimensional integral (not necessarily convex) polytopes in the subsection \ref{subsection_TMM}. To do this, we ``deform'' the Hermitian matrix model. Unfortunately the deformation violate the unitary invariance, but we can still apply some methods motivated by matrix models. From this point of view our models, strictly speaking, are not matrix (or tensor). But for the convenience, we will use this name in relation to the constructed models.

In the context of matrix models, it is natural to add all higher powers to the potential \cite{itoyama2017ward}. For example, considering Hermitian matrix model, we usually replace cubic potential in the following way:
\begin{equation}\label{renormalization}
    -\frac{1}{2}Tr(H^2)+t_3Tr(H^3)\rightarrow -\frac{\mu}{2} Tr(H^2)+\sum\limits_{n=0}^{\infty}t_nTr(H^n)
\end{equation}
We do the same for our matrix model of two-dimensional polytopes in subsection \ref{subsection_higher_degree}. It turns out that this is also natural from the point of view of integral polyhedra, since it enumerate subdivisions into polytopes. In addition, in subsection \ref{subsection_Virasoro} this allows us to write the Virasoro operator algebra of Ward identities in a similar way as for Hermitian matrix model (we briefly discuss the HMM example in subsection \ref{subsection_MM}):
\begin{equation}
    [L_n,L_m]=(n-m)L_{n+m}
\end{equation}
The presence of such an algebra of operators indicates the presence of integrability in such a model \cite{mironov1990origin,itoyama1991noncritical,ambjorn1990properties}.

Finally, in subsection \ref{subsection_higher_dimensions} we generalize this approach to the higher dimensional integral polyhedra, but in terms of ``tensor'' models \cite{freidel2005group,oriti2006group,boulatov1992model,akhmedov2005comment,guruau2017random}.

In section \ref{section_applications}, we present our considerations regarding the applications of the models we have obtained to the description of Batyrev's mirror pairs \cite{batyrev1994dual,batyrev1994hodge}. First, we give some basic definitions and statements from this area of science, after which we discuss what facts allow us to assume the applicability of our models to the enumeration of the Hodge numbers of Calabi-Yau hypersurfaces in toric varieties and the description of the moduli spaces of complex structures on them.

At the last section \ref{section_conclusion} we list our main results and discuss possible directions for further research on this topic and most interesting questions in our opinion.


\section{Preliminaries}\label{section_preliminaries}

In this section, we give some background information with references. Since our goal is to connect two subjects: matrix models and combinatorics of polyhedra, there are two subsections here, one for each side. We do not give any definitions and statements from toric geometry in this section, since in fact, the main discussion below is devoted to matrix models describing combinatorics of integral polyhedra. In section \ref{section_applications} we will discuss toric geometry and our thoughts on its connections with our matrix model.

\subsection{Convex integral polyhedra}

From the side of combinatorics, the object we are studying -- an integral convex polyhedron -- has a simple definition. Nevertheless, its combinatorics isn't trivial and quite deep. As mentioned in introduction, in future we would like to connect the theory of matrix models with the theory of generic hypersurfaces in toric varieties, so we interested in numbers $\ell(\Delta)$ and secondary polytope. We give both definitions here.

\vspace{0.5cm}

Let $M\simeq \mathbb{Z}^n$ be an $n$-dimensional lattice (free abelian group). The corresponding $\mathbb{R}$-vector space defined as $M_{\mathbb{R}}:=M\otimes \mathbb{R}\simeq \mathbb{R}^{n}$.

A \textit{convex polyhedron} in $M_{\mathbb{R}}$ is the convex hull of finitely many points in $M_{\mathbb{R}}$. A \textit{convex integral polyhedron} in $M_{\mathbb{R}}$ is the convex hull of finitely many points from $M$. The \textit{dimension} of a polytope is the dimension of the minimal affine subspace containing this polytope.

An affine hyperplane intersecting a polyhedron is called \textit{supporting} if the entire polyhedron lies in one of the two closed half-spaces defined by it. The intersection of such a plane with a polyhedron is called a \textit{face} (which, of course, is a convex polyhedron). If the dimension of a face is equal to $n-1$, then it is called \textit{facet}.

Let us define the numbers $\ell(\Delta)$ and $\ell^*(\Delta)$ (the notation is taken from Batyrev's papers on dual polyhedra and Calabi-Yau hypersurfaces) for a convex integral polytope $\Delta$. $\ell(\Delta)$ is just the number of integral points (i.e., $M$ lattice points) in the polyhedron:

\begin{equation}
    \ell(\Delta):=\#\{x\in M|x\in\Delta\}
\end{equation}

$\ell^*(\Delta)$ is the number of interior points of a polyhedron that do not belong to its facets.

As mentioned above, another combinatorial object that interests us is the secondary polytope \cite{gelfand1994discriminants}.

The easiest way to construct a secondary polytope is to use the characteristic functions of all possible $\Delta$ triangulations.

A \textit{triangulation} $\tau$ of the integral polytope $\Delta$ is a decomposition of it into a finite number of simplices (of dimension $n$), such that the vertices of the simplices are integral points, and the intersection of any two simplices is their common (possibly empty) face.

For each triangulation $\tau$ there is the associated \textit{characteristic function} $\phi_{\tau}$ on the set of vertices of this triangulation:

\begin{equation}
    \phi_{\tau}(a)=\sum\limits_{\{\sigma|a\in\sigma\}}\textrm{Vol}(\sigma),
\end{equation}

\noindent where the summation is over all simplices for which $a$ is a vertex.

Now one has to extend $\phi_T$ to all points in $\Delta\cap N$, taking $\phi_{\tau}(a)=0$ if $a\in\Delta\cap N$ is not a vertex of the triangulation $\tau$.

The \textit{secondary polytope} $\mathrm{Sec}(\Delta)$ of the polytope $\Delta$ can be defined as the convex hull of all vertices $(\phi_{\tau})$ in the vector space $\mathbb{R}^{\Delta\cap N}$:
\begin{equation}\label{secondary_polytope}
    \mathrm{Sec}(\Delta)=\textrm{Conv}(\{(\phi_{\tau})\in\mathbb{R}^{\Delta\cap N}|\tau\in \text{triangulations of } \Delta\}).
\end{equation}

Note here that the point $\phi_{\tau}$ is a vertex of the secondary polytope if and only if the triangulation is regular (or coherent, see \cite{gelfand1994discriminants}). The faces of the secondary polyhedron can be described in a similar way. They correspond to \textit{regular subdivisions} of a polytope into convex polygons. Below we construct the matrix model that enumerates subdivisions into polygons, which, however, may not be convex.

\subsection{Matrix models}\label{subsection_MM}

In this subsection we discuss some of the main concepts of matrix models, we are interested in, for the purposes of our paper. We demonstrate all of them by example of Hermitian matrix model (HMM) since our model constructed in subsection \ref{subsection_TMM} is a ``deformation'' of this model.

\vspace{0.5cm}

The partition function of HMM is the following:

\begin{equation}\label{HMM_partition_function}
    Z_N(t)=\int_{formal}\mathcal{D}H\exp\left[N\left(-\mu/2Tr(H^2)+\sum\limits_{q=0}^{\infty}t_qTr(H^k)\right)\right]
\end{equation}

\noindent where integration is over hermitian $N\times N$ matrices with the measure:
\begin{equation}\label{hermitian_measure}
    \mathcal{D}H:=\prod\limits_{i}dH_{i}^i\prod\limits_{i<j}dRe(H_{i}^j)dIm(H_{i}^j)
\end{equation}

The notation $\int_{formal}$ means that the exponent under the integral must first be expanded into a series in the variables $t=(t_0,t_1,\dots)$, and then each coefficient separately for each monomial of $t$ should be integrated.

We will be interested in expressions of the form (called correlators):

\begin{equation}\label{correlator}
    \langle\dots\rangle_0:=\frac{\int\mathcal{D}H\dots \exp[-N/2Tr(H^2)]}{\int\mathcal{D}H\exp[-N/2Tr(H^2)]}
\end{equation}

where $\dots$ is an arbitrary polynomial expression in matrix elements of $H$.

Since the measure in \eqref{correlator} is Gaussian, there is Wick's rule.

\vspace{0.5cm}

\textbf{Wick's rule.} The correlator of an arbitrary product of matrix elements $H_i^j$ is split into the sum:

\begin{equation}
    \langle H_{i_1}^{j_1}\dots H_{i_n}^{j_n}\rangle_0=\sum\limits_{\text{pairings}}\prod\limits_{\text{pairs}}\langle H_{i}^{j}H_{i'}^{j'}\rangle_0
\end{equation}
where the summation is carried out over all ways to split the factors into pairs, and the product is taken over all pairs. Furthermore:

\begin{equation}
    \langle H_{i}^{j}H_{k}^{l}\rangle_0=\frac{\delta_{i}^{l}\delta_{j}^{k}}{N}
\end{equation}

\textbf{Example.} Suppose we want to calculate the correlator $\langle Tr(H^4)\rangle_0$. Wick's rule allows one to do this in the following way:

$$\langle N Tr(H^4)\rangle_0=N(\langle H_i^jH_j^kH_k^lH_l^i \rangle_0=\langle H_i^jH_j^k\rangle_0\langle H_k^lH_l^i\rangle_0+\langle H_i^jH_k^l\rangle_0\langle H_j^kH_l^i\rangle_0+\langle H_i^jH_l^i\rangle_0\langle H_j^kH_k^l\rangle_0)=$$
$$=\frac{1}{N}(\delta^j_j\delta_i^k\delta^l_l\delta_k^i+\delta^j_k\delta^l_i\delta^k_l\delta^i_j+\delta^j_l\delta^i_i\delta^l_j\delta^k_k)=2N^2+1$$

All this makes it possible to give a geometric meaning to correlators by considering an analog of Feynman diagrams. In the case of Hermitian matrix model, one of the interpretations of the diagram technique has the following form. Omitting details, let us consider a correlator of the form

\begin{equation}
    \langle NTr(H^{k_1})\dots NTr(H^{k_n})\rangle_0
\end{equation}

Each factor $NTr(H^{k_i})$ is associated with a $k_i$-gon, and then all possible gluings of polygons along their edges into a closed orientable surface are considered. It turns out that each gluing contributes

\begin{equation}
    \sim N^{\chi(G)}
\end{equation}
where $\chi(G)$ is the Euler characteristic of the resulting surface.

\textbf{Example.} In the previous example, we assotiate a square to the expression under the correlator. The three terms below correspond to three ways of gluing the edges of the square in pairs. Two of them correspond to the sphere $S^2$ and contribute $2N^{\chi(S^2)}=2N^2$, while the last one corresponds to the torus $T$ and contributes $N^{\chi(T) }=1$.

Thus, the correlators of the form $\langle NTr(H^{2n})\rangle_0$ in the Hermitian matrix model enumerate the number of ways to glue a $2n$-gon along the edges into oriented surfaces of different genera.

\vspace{0.5cm}

Similar diagram techniques also occur in other matrix models. For example, in Kontsevich matrix model, summation is carried out over trivalent weighted graphs describing the cellular decomposition of the decorated moduli space of complex curves \cite{zvonkine2002strebel,Kontsevich_1992}.

The matrix model we are constructing uses a similar idea to enumerate triangulations of integral polyhedra.

\vspace{0.5cm}

\textbf{Ward identities.} One can write equations on correlators in matrix models that follow from the reparametrization invariance of the matrix integral. These equations are called Ward identities and often allow one to ``completely solve'' the matrix model, which means writing recursive relations on all correlators (invariant under the conjugation $H\rightarrow UHU^{+}$ by unitary matrix), allowing one to recover any correlator from a small a set of initial data.

Let's look at Hermitian matrix model. Any reparametrization of the form $H\rightarrow H+\delta H^n$ does not change the value of the integral, but transforms the integrand, which leads to the Ward identities. We will, however, use a different approach to derive them. It consists in the following. Substituting a differential operator of the form $\frac{\partial}{\partial H_{ij}}$ under the integral, we obtain 0, because the Gaussian exponent vanishes the integrand at infinity:

\begin{equation}
    0=\int\mathcal{D}H\frac{\partial}{\partial H_{i}^{j}}\dots \exp[-N/2 Tr(H^2)]
\end{equation}

\noindent where $\dots$ implies some polynomial expression of matrix elements of $H$.

Inserting the operator 
\begin{equation}\label{L_operator_partial}
    Tr\left(\frac{\partial}{\partial H}H^{n+1}\right)
\end{equation}
into the partition function \eqref{HMM_partition_function} gives the Ward identities, which can be written as:

\begin{equation}
    L_nZ=0, \:\geq -1,
\end{equation}

\noindent where

\begin{equation}\label{L_operator}
    L_n=-\mu\frac{\partial}{\partial t_{n+2}}+\sum\limits_{k=0}^{\infty}kt_k\frac{\partial}{\partial t_{k+n}}+\sum\limits_{k=0}^{n}\frac{\partial^2}{\partial t_k\partial t_{n-k}}+2N\frac{\partial}{\partial t_n}+N^2\delta_{n,0}
\end{equation}

Surprizing fact here is that such operators form the Virasoro algebra:

\begin{equation}
    [L_n,L_m]=(n-m)L_{n+m}
\end{equation}

The existence of such relations indicates the existence of a sufficient number of identities to recover arbitrary invariant correlators \cite{mironov1990origin}.

The similar operators and relations will appear in our model, but this will only allow us to recover the correlators, we are interested in, from the ``simple'' ones. A detailed discussion of this issue is provided in the section \ref{section_conclusion}.

\section{Matrix model of triangulations}\label{section_MM}

In this section, we build a matrix model, which we call the matrix model of triangulations. As the name suggests, the main idea is to consider all possible triangulations with vertices at integral points of integral polytope $\Delta$. We first consider two-dimensional polytopes and then generalize the model to higher dimensions. We think that it is quite natural to consider for our purposes all possible triangulations for the following reasons:

\begin{itemize}
    \item Matrix models with a potential of the 3rd degree are naturally enumerate triangulations of two-dimensional surfaces in their diagram technique \cite{KAZAKOV1985295}. Therefore, we only have to modify such a matrix model so that the vertices of the triangles correspond to integral points and the two-dimensional surface is simply a disk with a boundary coinciding with the boundary of the polytope.
    \item The number of points in a maximal triangulation of the polytope is expressed in terms of the number of triangles and doesn't depend on particular choice of such triangulation. Therefore, by introducing the parameter at the cubic potential in the matrix model, we can calculate the number of triangles, and, hence, the number of integral points $\ell(\Delta)$
    \item The secondary polytope \eqref{secondary_polytope}, which describes moduli spaces of toric hypersurfaces in some sense \cite{gelfand1994discriminants}, is the convex hull of the characteristic functions of all possible triangulations. Therefore, enumeration of all triangulations can be useful for studying moduli spaces
\end{itemize}

So we construct the simplex model in case of two dimensional polytopes in the subsection \ref{subsection_TMM}. Then we generalise it to higher degree potential in subsection \ref{subsection_higher_degree} which corresponds to subdivisions instead of triangulations. This generalization allows us to obtain Ward identities and operator \eqref{L_operator} analogue in the subsection \ref{subsection_Virasoro}. Finally, we generalize our model to higher dimensional case and do the same steps for it in the subsection \ref{subsection_higher_dimensions}.

\subsection{Constructing matrix model for polygons}\label{subsection_TMM}

In this subsection, the dimension of the lattice $M$ is equal two. It is well known that the matrix model with a potential of the third degree enumerates triangulations of two-dimensional (in the sense of real dimension) oriented surfaces \cite{KAZAKOV1985295}. We mentioned about it considering example of the Hermitian matrix model (subsection \ref{subsection_MM}). In our case, the following restrictions must be imposed on triangulations:
\begin{itemize}
    \item firstly, the triangulation vertices must belong to the two-dimensional lattice $M=\mathbb{Z}^{2}$
    \item secondly, the triangles must not overlap each other in the lattice $M=\mathbb{Z}^{2}$
    \item thirdly, the surface glued from triangles has the boundary $S^1$, which is the boundary of the polygon, and, furthermore, the surface has the topology of the disk, which is the interior of the polygon
\end{itemize}

\vspace{0.5cm}

Let us give some considerations on how we are going to construct the corresponding matrix model.

In order for the triangulation vertices to belong to the lattice, we define an ``index set'' $I^2_{N}:=\{m\in \mathbb{Z}^2:\: |m_1|,|m_2|\leq N \}$. This means that we will consider Hermitian matrices $(H_{ij})$ of size $(2N+1)^2\times (2N+1)^2$, whose indices $i$ and $j$ take values from the set $I ^2_N$. Hence the triangulation vertices will correspond to points from $I^2_N$. Taking the limit $N\rightarrow \infty$, we obtain the entire lattice $\mathbb{Z}^2$.

In order for triangulations to have boundaries we will consider correlators of special form:

\begin{equation}\label{polygon_correlator}
    \langle\Delta_{i_1,\dots,i_n}\rangle:=\left\langle H_{i_1i_2}H_{i_2i_3}\dots H_{i_ni_1}\right\rangle
\end{equation}

\noindent where the sequence $i_1,\dots,i_n\in I^2_N$ is formed by vertices of fixed polygon $\Delta_{i_1,\dots,i_n}$ in clockwise order (there is no summation over indices $i_1,\dots,i_n\in I^2_N$ here). Note that we can set not only a boundary of a polygon in this way and not only a boundary of a convex polygon. For now we will consider only boundaries of polygones, possibly not convex one. We leave the discussion of other correlators to further research.

Finally, in order for the triangles not to overlap, we require them to be oriented counterclockwise with respect to the lattice and to be non-degenerate. Then it is easy to check that if the boundary of the triangulation is given by the boundary of an integral polyhedron, then there will be no overlapping of triangles. Moreover, no triangles can arise outside the interior of a polygon. In order to take into account only properly oriented and non-degenerate triangles, we will have to significantly deform the cubic potential, writing it in the form 

\begin{equation}\label{W_potential}
    W(H):=\sum\limits_{i,j,k\in I^{2}_{N}}\omega_{ijk}H_{ij}H_{jk}H_{ki}
\end{equation}

\noindent where

\begin{equation}
    \omega^{ijk}\sim\begin{cases}
    1, \text{ for } \epsilon^{ab}(i_a-j_a)(k_b-j_b)>0 \\
    0, \text{ for } \epsilon^{ab}(i_a-j_a)(k_b-j_b)\leq 0
    \end{cases}.
\end{equation}
Recall that indices are takes values in the set $I^2_N$ which is a part of two-dimensional lattice, so they have two components $i=(i_1,i_2)$.

Thus, we will consider the following matrix model:

\begin{equation}\label{TMM}
    Z=\int_{formal} \mathcal{D}H \exp\left[-\frac{1}{2}Tr(H^2)+\beta W(H)\right].
\end{equation}

Here, $\int_{formal}$, as above, means that we consider the matrix model as a formal series in $\beta$, that is, the integral is taken for each degree of $\beta$ in the expansion of the exponent separately. Integration is over Hermitian matrices, in the sense that $H_{ij}=\overline{H}_{ji}$ and the measure $\mathcal{D}M$ has the same form as in \eqref{hermitian_measure}.

We intentionally write both indices of $H$ as subscripts, because the matrix model \eqref{TMM} is not unitary invariant. In this sense, our model is not formally a matrix model, we only use the motivation of matrix models to write the finite-dimensional integral in this form. On the other hand, it is possible to preserve the matrix form (and, in further subsections, the tensor structure) by writing similar expressions in a slightly more complex form, but this is not required for our purposes.

Note that $Z=0$ due to the orientation of all triangles and the absence of boundary in the sense of \eqref{polygon_correlator}.

In order to clarify now the connection between correlators \eqref{polygon_correlator} and integral polygons, we formulate the following theorem.

\vspace{0.5cm}

\textbf{Theorem 1.}\textit{ Let the potential $W(H)$ be defined by $\omega_{ijk}$ of the form }
\begin{equation}\label{omega}
    \omega_{ijk}=\begin{cases}
    e^{\epsilon^{ab}(i_a-j_a)(k_b-j_b)(x^i+x^j+x^k)}, \text{ for } \epsilon^{ab}(i_a-j_a)(k_b-j_b)>0 \\
    0, \text{ for } \epsilon^{ab}(i_a-j_a)(k_b-j_b)\leq 0
    \end{cases},
\end{equation}

\noindent\textit{where $x^i$ are formal variables ($i\in I^2_N$). Then the correlator} 

\begin{equation}\label{decomposition}
    \langle\Delta_{i_1,\dots,i_n}\rangle:=\left\langle H_{i_1i_2}H_{i_2i_3}\dots H_{i_ni_1}\right\rangle=\sum\limits_{k=0}^{\infty}\frac{\beta^k}{k!}\left\langle H_{i_1i_2}H_{i_2i_3}\dots H_{i_ni_1} (W(H))^k\right\rangle_0,
\end{equation}

\noindent\textit{where the sequence of points $i_1,\dots,i_n\in I^2_N$ forms the boundary of the polygon $\Delta_{i_1,\dots,i_n}$, for sufficiently large $N$ is expressed in terms of the sum over triangulations $\Delta_{i_1,\ dots,i_n}$ of the following form:}
\begin{equation}
    \langle\Delta_{i_1,\dots,i_n}\rangle=\sum\limits_{\tau\in \textrm{triangulations of }\Delta}\beta^{|\tau|}e^{2\langle\phi_{tau},L(\Delta)\rangle}
\end{equation}
\textit{where $|\tau|$ is the number of triangles in the triangulation $\tau$, $\langle\phi_{tau},L(\Delta)\rangle:=\sum\limits_{a\in V(\tau)}\phi_{a}(\tau)x^a$ is a polynomial whose Newton polyhedron is $\Delta$, and whose coefficients are the values of the triangulation $\tau$ characteristic function at the corresponding triangulation vertex $a\in V(\tau) $.}

Let us immediately note an important corollary which shows how the numbers $\ell(\Delta)$ arise in this model.

\vspace{0.5cm}

\textbf{Corollary.} \textit{The highest degree of $\beta$ in the correlator $\langle\Delta_{i_1,\dots,i_n}\rangle$ is $2\ell^*(\Delta)+n-2$.}

\subsection{Higher degree generalization}\label{subsection_higher_degree}

There is some motivation to generalize the potential of the matrix model to higher degrees (see \cite{itoyama2017ward}). Adding higher degree terms to potential can be interpreted as ``dressing process'' in terms of renormalization group theory. Such a procedure naturally leads to a Virasoro-like algebra of operators corresponding to the Ward identities. Here we present a way to generalize the above model in a similar way: we introduce higher powers into the matrix model \eqref{TMM}. We also discuss how the Theorem 1 changes. 

\vspace{0.5cm}

Let us consider the following matrix model (inspired by analogous replacement \eqref{renormalization} for HMM):

\begin{equation}\label{MM_subdivisions}
    Z=\int \mathcal{D}H \exp\left\{-\frac{\mu}{2} Tr(H^2)+\sum\limits_{n=1}^{\infty}t_n W_{n}(H)\right\}
\end{equation}

\noindent where

\begin{equation}
    W_{n}(H)=\sum\limits_{i_1,\dots,i_{n+2}\in I^2_N}\omega^n_{i_1\dots i_{n+2}}H_{i_1i_2}H_{i_2i_3}\dots H_{i_{n+2}i_1}
\end{equation}

We define $\omega^n_{i_1\dots i_{n+2}}$ by recursion relations having sense of gluing triangles corresponding to $\omega^1_{ijk}$ into polygons:

\begin{itemize}
    \item $\omega^1_{ijk}$ is just $\omega_{ijk}$ defined in previous subsection. It correspond to the triangle with vertices $i,j,k\in I^2_N$
    \item Cyclic permutation of indices defines the same $n+2$-gon
    \begin{equation}
        \omega^{n-2}_{i_1\dots i_{n-1}i_{n}}=\omega^{n-2}_{i_ni_1\dots i_{n-1}}
    \end{equation}
    \item Gluing two polygons along a common edge gives a new polygon
    \begin{equation}
        \omega^{m+n-4}_{i_1\dots i_nj_3\dots j_m}=\begin{cases}
            \omega^{n-2}_{i_1\dots i_n}\omega^{m-2}_{i_1i_nj_3\dots j_m},\text{ if } \Delta_{i_1\dots i_nj_3\dots j_m} \text{ is convex}\\
            0, \: \text{ otherwise}
        \end{cases}
    \end{equation}
\end{itemize}

We can generalize Theorem 1 for this matrix model as follows:

\vspace{0.5cm}

\textbf{Theorem 1'.} \textit{ Let the potential $W_1(H)$ be of the form \eqref{W_potential}, and the remaining $W_n(H)$ satisfy the recursion relations described above. Then the correlator} 

\begin{equation}
    \langle\Delta_{i_1,\dots,i_n}\rangle:=\left\langle H_{i_1i_2}H_{i_2i_3}\dots H_{i_ni_1}\right\rangle,
\end{equation}

\noindent\textit{where the sequence of points $i_1,\dots,i_n\in I^2_N$ forms the boundary of the polygon $\Delta_{i_1,\dots,i_n}$, for sufficiently large $N$ is expressed in terms of the sum over subdivisions into convex polygons of the polygon $\Delta_{i_1,\dots,i_n}$ of the following form:}
\begin{equation}\label{main_th}
    \langle\Delta_{i_1,\dots,i_n}\rangle=\sum\limits_{ \textrm{subdivisions of }\Delta}\#_{\text{subdivision}}1/\mu^{\# \text{ of subdivision edges}}\prod\limits_{q=1}^{\infty} t_q^{\# \text{ of q-gons in subdivision}}\prod\limits_{\square\in\text{ q-gons}}W_q(\square)
\end{equation}
\textit{where $W_q(\square)$ is the value of the potential $W_q(H)$ on a particular polygon $\square$. $\#_{\text{subdivision}}$ is some coefficient depending on subdivision.}

\subsection{Ward identities/Virasoro algebra}\label{subsection_Virasoro}

Here we show how to obtain the Virasoro constraints for the matrix model \eqref{MM_subdivisions}.

\vspace{0.5cm}

We deform operators \eqref{L_operator_partial} as follows. From the Feynman diagram point of view this operator corresponds to insertion of an $n$-gon:

\begin{equation}\label{deformed_L_operator}
    \hat{L}_{n}=\sum\limits_{i_1\dots i_{n+2}}\omega^{n}_{i_1i_2\dots i_{n+2}}\frac{\partial}{\partial M_{i_2 i_1}}M_{i_2i_3}M_{i_3 i_4}\dots M_{i_{n+2} i_1}
\end{equation}

As it was for Hermitian matrix model we obtain Ward identities:

\begin{equation}\label{Virasoro_constraints}
    0=\hat{L}_{n-1}\langle\Delta_{j_1\dots j_m}\rangle
\end{equation}

We act by operators $\hat{L}_{n-1}$ on a correlator of the form \eqref{polygon_correlator}, because it can not be recovered from the partition function by the action of any differential operator in variables $t_i$.

\vspace{0.5cm}

\textbf{Theorem 2.} \textit{The action of the operator \eqref{deformed_L_operator} can be written in the following form:}
\begin{equation}
    \hat{L}_{n}\langle\Delta_{j_1\dots j_m}\rangle=\left(-\mu\frac{\partial}{\partial t_{n+2}}+\sum\limits_{k=1}^{\infty}kt_k\frac{\partial}{\partial t_{k+n+2}}+\hat{I}_n\right)\langle\Delta_{j_1\dots j_m}\rangle
\end{equation}

\noindent\textit{where}

\begin{equation}
    \hat{I}_n\langle\Delta_{j_1\dots j_m}\rangle=\sum\limits_{l}\sum\limits_{i_1,\dots ,i_{n}}\omega^n_{j_li_1\dots i_{n}j_{l+1}}\langle M_{j_1j_2}\dots \hat{M}_{j_lj_{l+1}}\dots M_{j_mj_1}(M_{j_li_1}M_{i_1i_2}\dots M_{i_nj_{l+1}})\rangle=
\end{equation}
$$=\sum\limits_{l}\sum\limits_{i_1,\dots ,i_{n}}\omega^n_{j_li_1\dots i_{n}j_{l+1}}\langle \Delta_{j_1\dots j_li_1\dots i_n j_{l+1}\dots j_m}\rangle$$

\textit{And operators $\hat{L}_{n}$ form the Virasoro algebra:}

\begin{equation}\label{Virasoro_algebra}
    [\hat{L}_{n},\hat{L}_{m}]=(n-m)\hat{L}_{n+m}
\end{equation}

Note here that correlators obtained by the action of operator $\hat{I}_n$ does not necessarily correspond to boundaries of polygons. So we have to consider correlators of the form \eqref{polygon_correlator} for an arbitrary sequence of points $i_1,\dots,i_n\in I^2_N$ in order to obtain recursion relations from Virasoro constraints \eqref{Virasoro_constraints}.

\subsection{Higher dimensional generalization and tensor model}\label{subsection_higher_dimensions}

In this subsection, we give some considerations about higher dimensional generalizations of the matrix model \eqref{TMM}. We are going to construct a ``matrix model'' that enumerates subdivisions of $d$-dimensional polytopes. But in order to describe the gluing of simplices of higher dimensions, one has to write a tensor model instead of a matrix one. Such models are poorly understood, and in general it is not known whether their Ward identities are solvable. However, in our case, this is not very significant, since as in the model described above, the tensor structure would not be important to us. Indeed, many methods of working with matrix models are based on the possibility to diagonalize the matrix and integrate it over eigenvalues or, for example, use the Itzykson-Zuber formula. Both of these methods are based on the unitary invariance of the action, which is missing in our model and is not used in the previous reasoning. Therefore, we expect that all previous statements can be generalized to the case of a tensor model.

Note that there are some multidimensional generalizations of matrix models which describe higher dimensional triangulations of pseudo manifolds \cite{gurau2010lost}. Such generalizations are called group field theores (GFT) \cite{freidel2005group,oriti2006group,guruau2017random,tanasa2021combinatorial}. One can consider our tensor model as a deformation of GFT in the same sense that the matrix model of triangulations is the deformation of Hermitian one. The Virasoro-like constraints \eqref{Virasoro_like} we obtain are quite similar to the Virasoro-like constraints which arise in the large $N$ limit for GFT. They are also can be naturally described in terms of trees \cite{gurau2011generalization}. Nonetheless, there is a difference in how we glue trees and how they are glued in GFT. It would be interesting to understand how these two Virasoro-like constraints are connected, but we leave that for further research.

\vspace{0.5cm}

Now we fix the set in which the indices of tensors take values, in the following obvious way $I_N^d=\{(m_1,\dots,m_d)\in \mathbb{Z}^d:|m_i|\leq N\:\forall i\}$.

We associate the $d-1$-simplex with vertices $i_1\dots i_d$ with the tensor $T_{i_1 \dots i_d}$ as it was for our matrix model: $H_{ij}$ corresponds to the edge of the triangle with vertices $ij$. In addition, we impose the following condition on the tensor $T_{i_1\dots i_d}$:

\begin{equation}
    T_{i_1\dots i_li_{l+1}\dots i_d}=\overline{T}_{i_1\dots i_{l+1}i_l\dots i_d}, \:\forall l\in \{1,\dots,d-1\}
\end{equation}

Then the analogue of the quadratic term will be:
\begin{equation}
    -1/2\sum\limits_{i_1\dots i_n} T_{i_1\dots i_{d-1} i_d} T_{i_2i_1i_3\dots i_d}
\end{equation}

In order to write an analogue of the cubic action, we have to ``glue'' the $d$-dimensional simplex from its boundaries. For example, the cubic potential associated with 3-dimensional simplex can be written in the following way:
\begin{equation}
    \sum\limits_{i_1i_2i_3i_4}\Omega^{4}_{i_1i_2i_3i_4} T_{i_2 i_3 i_4}T_{i_1 i_3 i_4}T_{i_1 i_2 i_4}T_{i_1 i_2 i_3}
\end{equation}

In general, this potential can be written as:

\begin{equation}\label{qubic_term}
    \sum\limits_{i_0\dots i_{d}}\Omega^{d+1}_{i_0\dots i_d}\prod\limits_{k=0}^{d}T_{i_0\dots\hat{i}_k\dots i_d}
\end{equation}

If we want to obtain an analogue of the Ward identities \eqref{Virasoro_constraints} and the Virasoro algebra \eqref{Virasoro_algebra}, then we have to add higher powers in $T$ to the potential. However, things are more complicated for them then it was for the matrix model. We argue that they must correspond to ``tree'' simplicial complexes instead polygons arising in the matrix model. We mean the following.

\begin{itemize}
    \item This simplicial complex is a simplicial manifold of dimension $d$ with boundary, i.e. it is glued from $d$-simplices so that for each $d-1$-simplex, at most 2 $d$-simplices are glued (in other words it is a triangulation of $d$-dimensional manifold with triangulated boundary). The boundary of such a simplicial complex is a simplicial sphere.
    \item Therefore,  the dual simplicial complex can be associated with it. Moreover, this dual simplicial complex has to be a tree with ordered edges at each vertex.
\end{itemize}

Therefore, for each degree of the potential, we have to write the sum over all possible trees of this form with a fixed number of vertices. In the case of the matrix model, this problem was greatly simplified, since all such trees give equivalent partitions of the sphere $S^1$ into simplices. It was expressed in the possibility of a cyclic permutation of $\omega^n_{i_1\dots i_{n+2}}$ indices.

Let $\mathcal{T}$ be the set of trees described above. For each tree $\tau\in \mathcal{T}$ there is the associated polytope $\Delta(\tau)$ and associated term of potential $W_{\tau}(T)$. The partition function of our tensor model has the following form:

\begin{equation}\label{Tensor_model}
    Z(t)=\int\mathcal{D}T \exp\left\{-1/2\sum\limits_{i_1\dots i_n} T_{i_1\dots i_{d-1} i_d} T_{i_2i_1i_3\dots i_d}+\sum\limits_{\tau\in\mathcal{T}}t_{\tau}W_{\tau}(T)\right\}
\end{equation}

\noindent with the set of variables $\{t_{\tau}|\tau\in\mathcal{T}\}$.
We give a more detailed description of the tensor model potential, trees and associated $\hat{L}_{(\tau,v)}$ operators with them in Appendix B. Here we only present the final results.

\textbf{Proposition 1.} (Generalization of the Theorem 1'). \textit{Let $\Delta$ be integral simplicial polytope of dimension $d$ with the set of facets $F(\Delta)$. And let $V(\sigma)\subset I^d_N,\:\sigma\in F(\Delta)$ denote the set of vertices of the simplex $\sigma$. We associate the correlator}

\begin{equation}
    \langle\Delta\rangle = \left\langle\prod\limits_{\sigma\in F(\Delta)}T_{V(\sigma)}\right\rangle
\end{equation}

\noindent\textit{to the polytope $\Delta$, where $T_{V(\sigma)}:=V_{i_1\dots i_d},\:V(\sigma)=\{i_1,\dots,i_d\}$. Then (for sufficiently large $N$) this correlator can be expressed in terms of sudivisions of $\Delta$ in the following way:}

\begin{equation}
    \langle \Delta\rangle = \sum\limits_{\text{subdivisions' of } \Delta } \#_{\text{subdivision}}1/\mu^{\# \text{ of } d-1\text{-simplices}}\prod\limits_{\tau\in \mathcal{T}}t_{\tau}^{\#\Delta(\tau)}W_{\tau}(\Delta(\tau))
\end{equation}

\noindent\textit{where summation is over all subdivisions into convex polytopes associated with trees $\mathcal{T}$. $\#\Delta(\tau)$ is the number of polytopes $\Delta(\tau)$ in corresponding subdivision. $W_{\tau}(\Delta(\tau))$ is the value of  $W_{\tau}(T)$ on the polytope $\Delta(\tau)$. $\#_{\text{subdivision}}$ is some coefficient depending on subdivision.}

Unfortunately, we obtain here generating function only for subdivisions of special form. There exist polytopes of dimension $\geq 3$ that can not be represented as $\Delta(\tau)$ for some $\tau\in \mathcal{T}$. The simplest example is octahedron.

\textbf{Proposition 2.} (Generalization of the Theorem 2) \textit{The $\hat{L}_{(S,\sigma)}$ operators generate the Ward identities:}

\begin{equation}\label{tensor_Ward_identities}
    \hat{L}_{(S,\sigma)}\langle\Delta\rangle =0
\end{equation}

\textit{and form a Virasoro-like algebra}

\begin{equation}\label{Virasoro_like}
    [\hat{L}_{(\tau_1,v_1)},\hat{L}_{(\tau_2,v_2)}]=\sum\limits_{\tau\in(\tau_1\#(\tau_2,v_2)|v_1)}\hat{L}_{(\tau,v_1)}-\sum\limits_{\tau\in(\tau_2\#(\tau_1,v_1)|v_2)}\hat{L}_{(\tau,v_2)}
\end{equation}

\noindent\textit{where $(\tau_1\#(\tau_2,v_2)|v_1)$ is the set of all trees obtained by gluing the tree $\tau_2$ along the edge $v_2$ to an arbitrary edge of $\tau_1$ except $v_1$; the marked edge in resulting tree is $v_1$}

\section{Toric geometry applications}\label{section_applications}

In this section, we give some basic information about the connection between combinatorics of integral polyhedra and toric geometry. We hope this will clarify our motivation to study particular combinatorial objects described above.

\subsection{Toric geometry}

Toric varieties correspond one-to-one with the so-called fans \cite{fulton1993introduction}. That is why a lot of properties of toric varieties can be reformulated in combinatorial terms of fans. Let's give some examples.
\begin{itemize}
    \item A variety is complete (analytically compact) if and only if the union of all cones in the fan coincides with the ambient lattice.
    \item Toric variety is non-singular if and only if each cone of the corresponding fan is generated by a part of the basis (that is, any integral point of the cone can be obtained by a linear combination of generators that are part of the basis of the lattice with positive integer coefficients).
    \item Each cone of dimension $l$ corresponds to an orbit of a torus of dimension $n-l$ and an invariant subvariety of dimension $n-l$, which are the closures of these orbits.
    \item One-dimensional cones correspond to $T$-invariant divisors.
    \item Piecewise polynomial functions on the fan $\Sigma$ form a ring isomorphic to the equivariant Chow ring $A^{*}_T(X_{\Sigma})$ of the corresponding toric variety \cite{payne2006equivariant}.
\end{itemize}

We are interested in complete projective toric varieties. It is more convenient to describe such varieties in terms of convex integral polyhedra. For each polytope of this kind, one can construct a fan and a strictly convex piecewise linear function on it, which will define, respectively, a toric variety and an ample divisor on it (or, which is the same, an embedding in the projective space $\mathbb{C}P^k$ of sufficiently large dimension). In fact, a toric variety is projective if and only if the fan comes from a polytope (generally it is not true, the first counterexample occurs in dimension 4). That is why the main combinatorial object studied in this paper is an integral polyhedron.

Let us give here one of the ways to construct a projective toric variety from the polytope $\Delta$. Consider the ring $\mathbb{C}[x_0,x_1^{\pm 1},\dots,x_n^{\pm 1}]$, and the subring $S$ generated by monomials of the form $x_0^{m_0}x_1^ {m_1}\dots x_n^{m_n}$, where $(\frac{m_1}{m_0},\dots,\frac{m_n}{m_0})\in\Delta$. Then the corresponding toric variety
\begin{equation}\label{Proj_construction}
    X_{\Delta}=\textrm{Proj}(S).
\end{equation}

An interesting application of the theory of toric varieties is Batyrev's construction of Calabi-Yau mirror varieties. We give here some considerations from his work.

Let $\Delta\subset M_{\mathbb{R}}$ be a convex polytope, then the dual one defined as:

\begin{equation}
    \Delta^{\vee}=\left\{x\in M_{\mathbb{R}}|\langle x,y\rangle\geq -1,\:\forall y\in\Delta\right\}.
\end{equation}

Obviously, $((\Delta)^{\vee})^{\vee}=\Delta$.

The combinatorics of the dual polytope $\Delta^{\vee}$ is simply related to the combinatorics of the initial one $\Delta$. Each face $\Theta$ of dimension $k$ in $\Delta$ corresponds to a dual face $\hat{\Theta}$ in $\Delta^{\vee}$ of dimension $n-1-k$. And the face embedding relation is inverted: if $\Theta\subset \Gamma$ for $\Delta$, then $\hat{\Gamma}\subset \hat{\Theta}$ for $\Delta^{\vee}$.

In the context of constructing Calabi-Yau varieties, reflexive polyhedra are of particular interest. An integral convex polytope $\Delta$ is called reflexive if $\Delta^{\vee}$ is integral.

Calabi-Yau varieties arise as hypersurfaces in toric projective varieties $X_{\Delta}$, where $\Delta$ is reflexive. Strictly speaking, these hypersurfaces are given as divisors of global sections of the anticanonical sheaf in some desingularization of the toric variety associated with the reflexive polytope. 

However, there is a more explicit way to describe these hypersurfaces. Associate each node $m$ of lattice $M$ with coordinates $(m_1,\dots,m_n)$ with the Laurent monomial \\$x_1^{m_1}x_2^{m_2}\dots x_n^{m_n}$ and denote it by $ x^m$. The support of the Laurent polynomial $L=\sum\limits_{m\in M}a_m x^m$ with coefficients $a_m\in\mathbb{C}$ is a subset of the lattice $M$ of the form:

\begin{equation}
    \mathrm{Supp}(L)=\{m\in M|a_m\not=0\}.
\end{equation}

The Newton polytope of the polynomial $L$ is the convex hull of the support. The Newton polytope is obviously an integral polytope in $M_{\mathbb{R}}$.

Denote by $\mathbf{L}(\Delta)$ the space of polynomials with a fixed Newton polyhedron $\Delta$. Any polynomial $L$ from $\mathbf{L}(\Delta)$ is associated with the hypersurface $Z_{L,\Delta}=\{X\in T|L(X)=0\}$ in the torus and its closure $ \overline{Z}_{f,\Delta}$ in a projective toric variety $X_{\Delta}$. The closure $\overline{Z}_{f,\Delta}$ in the maximal desingularization of the variety $X_{\Delta}$ is exactly the Calabi-Yau hypersurface.

\vspace{0.5cm}

One of the consequences of mirror symmetry \cite{Witten1991zz} is that for a mirror pair $X$ and $X^{\circ}$ of Calabi-Yau varieties, the following Hodge number equality holds:

\begin{equation}\label{hodge_numbers}
    h^{1,1}(V)=h^{dim(V^{\circ})-1,1}(V^{\circ}), \:  h^{1,1}(V^{\circ})=h^{dim(V)-1,1}(V)
\end{equation}

Batyrev proved in his papers that pairs of Calabi-Yau hypersurfaces constructed from reflexive polyhedra are candidates for a mirror pair in the sense that the equalities \eqref{hodge_numbers} hold. It follows from explicit combinatorial formulas for the Hodge numbers of hypersurfaces:

\begin{equation}\label{hodge_11}
    h^{1,1}(V)=\ell(\Delta^{\vee})-n-1-\sum\limits_{\Gamma^{\vee}}\ell^*(\Gamma^{\vee})+\sum\limits_{\Theta^{\vee}}\ell^*(\Theta^{\vee})\ell^*(\hat{\Theta}^{\vee}),
\end{equation}
\begin{equation}\label{hodge_n1}
    h^{n-2,1}(V)=\ell(\Delta)-n-1-\sum\limits_{\Gamma}\ell^*(\Gamma)+\sum\limits_{\Theta}\ell^*(\Theta)\ell^*(\hat{\Theta}),
\end{equation}

\noindent where $\ell(\Delta)$ is the number of integral points in $\Delta$ mentioned above; $\Gamma$ and $\Gamma^{\vee}$ are faces of codimension one of the polytope $\Delta$ and its dual $\Delta^{\vee}$, respectively; $\Theta$ and $\Theta^{\vee}$ are faces of codimension two in $\Delta$ and $\Delta^{\vee}$, respectively; $\hat{\Theta}$ is a face of dimension one in $\Delta^{\vee}$ corresponding to the face $\Theta$ of codimension two in $\Delta$.

Note that some formulations of mirror symmetry expect the Hodge diamond of one variety to be the reflected Hodge diamond of another. As far as we know, the calculation of all Hodge numbers of Calabi-Yau hypersurfaces in Batyrev's construction is still an open question.

\subsection{Moduli space}\label{subsection_moduli}

A more general view on mirror symmetry is that there must be an isomorphism of moduli spaces of complex and Kähler structures of Calabi-Yau mirror varieties (with fixed Kähler and complex structures, respectively). Therefore, it would be interesting for us to touch upon the combinatorics associated with these moduli spaces. Here we give some known facts about the structure of moduli spaces of complex structures on Calabi-Yau hypersurfaces in toric varieties. In fact, we will consider a simplified moduli space of hypersurfaces in a toric variety $\overline{\mathcal{M}}_{simp}$. It is naturally described in terms of a secondary polyhedron \cite{gelfand1994discriminants}, on the one hand, and partially describes the real moduli space of complex structure on the Calabi-Yau hypersurface, on the other. Nevertheless, the moduli space of hypersurfaces is interesting in itself.

\vspace{0.5cm}

Let us give a brief description of connection between the moduli space of complex structures and $\overline{\mathcal{M}}_{simp}$, moving from the first to the second. First, to begin with, it is simpler to consider only the moduli space of complex structures arising from the embedding of a hypersurface in a toric variety denoted by $\mathcal{M}_{poly}$. In other words, it is the closure of the space of hypersurfaces in general position in a toric variety up to isomorphisms. It is natural to try to obtain this space from the space of Laurent polynomials $\mathbf{L}(\Delta)$. Recall that each (general position) polynomial in $\mathbf{L}(\Delta)$ corresponds to the hypersurface $Z_{f,\Delta}$ in the torus $T$ and its closure $\overline{Z}_{f, \Delta}$ in the toric variety $X_{\Delta}$.

One can describe the whole group of automorphisms of the toric variety $\mathrm{Aut}(X_{\Delta})$, whose action extends to the projective space $\mathbb{P}(\mathbf{L}(\Delta))$. This allows to define the following space\cite{batyrev1994hodge}:

\begin{equation}
    \mathcal{M}_{poly}:=\mathbb{P}(\mathbf{L}(\Delta))/\mathrm{Aut}(X_{\Delta}).
\end{equation}

It can be proved that the spaces $\mathcal{M}$ and $\mathcal{M}_{poly}$ coincide if the faces of $\Delta$ of codimension 2 do not contain integral points in their interior. Otherwise $\mathcal{M}_{poly}\subset\mathcal{M}$. This is consistent with the \eqref{hodge_n1} formula in the following sense:

\begin{itemize}
    \item $\textrm{dim}(\mathcal{M}_{poly})=\ell(\Delta)-n-1-\sum\limits_{\Gamma}\ell^*(\Gamma)$,
    \item $\textrm{dim}(\mathcal{M})=h^{n-2,1}=h^{n-2,1}(V)=\ell(\Delta)-n-1-\sum\limits_{\Gamma}\ell^*(\Gamma)+\sum\limits_{\Theta}\ell^*(\Theta)\ell^*(\hat{\Theta})$,
\end{itemize}

That is, the dimensions of the moduli spaces $\mathcal{M}$ and $\mathcal{M}_{poly}$ differ by a term that vanishes exactly when the set of interior integral points of codimension 2 faces is empty.

An even simpler object is the simple moduli space $\mathcal{M}_{simp}$ \cite{aspinwall1993monomial}, which differs from $\mathcal{M}_{poly}$ in that not the entire automorphism group is considered, but only its toric subgroup.

\begin{equation}
    \mathcal{M}_{simp}=\mathbb{P}(\mathbf{L}(\Delta)_0)/T,
\end{equation}

\noindent where $\mathbf{L}(\Delta)$ is the set of Laurent polynomials $L$ with support $\mathrm{Supp}(L)\subset \Delta\cap N$ such that the interior integral points of the facets of $\Delta$ do not lie in $\mathrm{Supp}(L)$.

It is known \cite{aspinwall1993monomial} that $\mathcal{M}_{simp}$ is a finite covering of $\mathcal{M}_{poly}$. For us, the following statement \cite{gelfand1994discriminants} will be important:

\vspace{0.5cm}

\textit{The compactification of the simple moduli space $\overline{\mathcal{M}}_{simp}$ is the toric variety associated with the secondary polytope of the polytope $\Delta$.}

\vspace{0.5cm}

This means that the simple moduli space can again be described in combinatorial terms. Since each face of the polyhedron corresponds to a closed orbit of the torus in the toric variety, and the set of such orbits forms a stratification of the toric variety, then the set of all regular subdivisions of the polytope $\Delta$ will correspond to the stratification of the moduli space $\overline{\mathcal{M}}_{simp}$. Therefore, it would be extremely interesting to modify the matrix model \eqref{MM_subdivisions} or tensor model \eqref{Tensor_model} in such a way that they enumerate precisely regular convex subdivisions. In this sense, it would resemble Kontsevich's matrix model, which enumerates in a certain sense the cellular structure of the decorated moduli space of complex structures.

\section{Discussion}\label{section_conclusion}

\subsection{Main results}

We constructed the matrix model of triangulations \eqref{TMM} for two-dimensional integral polyhedra. Correlators of the special form \eqref{polygon_correlator} in this model are generating functions for the numbers $\ell^*(\Delta)$ and the characteristic functions of triangulations in the sense of formula \eqref{main_th}. Then a generalized model \eqref{MM_subdivisions} is constructed, which enumerates already arbitrary subdivisions of integral polygones in a similar way (Theorem 1'). Due to the special structure of the potential in this model, a certain family of Ward identities can be written as operators \eqref{deformed_L_operator}, which form the Virasoro algebra.

Next, we gave considerations about a natural generalization to the case of integral polytopes of an arbitrary dimension. In this case, we have to consider the tensor models described in subsection \ref{subsection_higher_dimensions}. Despite the fact that we do not give an explicit form of the potential in such a model, it can be defined recursively (Appendix B). Generalization of the Theorem 1' is given: correlators of the tensor model enumerate subdivisions of polytopes in the same way as it the matrix model so. But only subdivisions of special type arise here (Proposition 1). Finally, the Ward identities are generalized to higher dimensional case, which in turn no longer form the Virasoro algebra, but lead to similar relations. Such a deformation of the Virasoro algebra seems rather interesting, but we leave it for further research.

\subsection{Open questions and further research}

The matrix models approach for studying the combinatorics of integral polytopes and the connection of this combinatorics with the geometry of toric varieties and hypersurfaces allows us to formulate many interesting questions and directions for further research. We list here an incomplete list of some of the most amusing questions in our opinion:

\begin{itemize}
    \item As mentioned above, the existence of the Virasoro algebra indicates the existence of ``complete solvability'' of the matrix model. So the next natural step would be to obtain explicit recurrent relations on the correlators of the matrix model \eqref{MM_subdivisions} for two-dimensional polytopes and, possibly, for the tensor model \eqref{Tensor_model}. We do not think, however, that in our case it will be possible to reconstruct an arbitrary correlator from a small dataset. The fact is that the ``complete solvability'' of a matrix model (for example, the Hermitian  one) means the possibility to calculate an arbitrary invariant (with respect to conjugation by unitary matrices) correlator. In our case, there is no unitary invariance, and we are interested in rather special correlators \eqref{polygon_correlator}. We assume that the recurrence relations will allow us to recover the generating function \eqref{main_th} for an arbitrary polyhedron from the same generating function for an arbitrary simplex.However, calculating such a generating function even for a simplex is a difficult task.
    \item We have already said that our models, strictly speaking, are not matrix or tensor, they do not have unitary invariance. However, this can be fixed by adding more tensor ``fields''. The question is whether these fields can be given a natural combinatorial/geometric meaning in terms of the diagram technique. It is also not clear what happens to the Ward identities in this case.
    \item Recall that our global goal is to describe the Hodge numbers and the moduli space of toric hypersurfaces using our matrix/tensor model. Despite the fact that the Hodge numbers are expressed in terms of the numbers $\ell(\Delta)$ \eqref{hodge_n1}, which are reproduced in our matrix model, we would like to obtain operators that, when acting on correlators \eqref{polygon_correlator}, would give exactly the Hodge numbers.
    \item As for the moduli spaces, first we would like to describe the stratification of the space $\overline{\mathcal{M}}_{simp}$. For this, as mentioned above, it is necessary to consider regular subdivisions into convex polyhedra. Therefore, we would like to see if it is possible to modify our matrix/tensor models so that they reproduce only such subdivisions.
    \item Another question is how to describe reflexive pairs of polytopes in terms of matrix models. Is it possible to write an operator that, when acting on correlators \eqref{polygon_correlator}, does not vanish only the correlators of reflexive polyhedra? Maybe such a result can be obtained in a certain limit of our models?
    \item 
    We do not know how to describe trees which gives the same polytope (or corresponding term in potential). It would be interesting to understand whether there exist such elementary tree transformations that preserve associated polytopes and generate all of such transformations.
\end{itemize}

\section{Acknowledgements}
Our work is supported in part by the grant of the Foundation for the Advancement of Theoretical Physics “BASIS” and by RFBR grants 20-01-00644 and 21-51-46010. We are grateful to Alexey Sleptsov and Alexandr Popolitov for very useful discussions and remarks.

\newpage

\section{Appendix A. Proof of the Theorem 1}

\textit{Proof of the Theorem 1 and it's Corollary}

As mentioned in subsection \ref{subsection_TMM}, the diagram technique corresponds to gluing triangles with vertices at integral points of the bounded lattice $I^2_N$. Indeed, to each factor $\omega^{ijk}H_{ij}H_{jk}H_{ki}$ under $\langle\dots\rangle_0$ we associate a trivalent vertex of the ribbon graph with weight $\omega^{ijk}$ , whose three boundaries correspond to indices $i,j,k$. The multiplier $H_{i_1i_2}\dots H_{i_ni_1}$ is associated with a vertex of valency $n$ with indices $i_1,\dots,i_n$ on the boundaries. Averaging with the Gaussian weight $\langle\dots\rangle_0$ according to Wick's rule yields all possible gluings of such vertices into an orientable ribbon graph. After that, we can continue the ribbon graph to a two-dimensional oriented surface by gluing two-dimensional disks along all boundaries of it. Therefore, we obtain the sum over the genus $g$ of the surface and all graphs $\Gamma(g)$ on them with one vertex of valence $n$ and $k$ vertices of valence 3 (where $k$ is the degree of $\beta$):

\begin{equation}\label{Wick_decomposition}
    \left\langle H_{i_1i_2}H_{i_2i_3}\dots H_{i_ni_1} (W(H))^k\right\rangle_0=\sum\limits_{g=0}^{\infty}\sum\limits_{\gamma\in \Gamma(g)}\sum\limits_{\{i_f\in I^2_N,\:f\in F(\gamma)\}}\prod\limits_{v\in V'(\gamma)}\omega_v^{i_{f_1(v)}i_{f_2(v)}i_{f_3(v)}}
\end{equation}
\noindent where $F(\gamma)$ is the set of faces of the graph $\gamma$ on the corresponding surface, $V'(\gamma)$ is the set of vertices, except for the vertex of valency $n$, $f_{j}(v)$ are graph faces on whose boundary lies the vertex $v$ ($j\in\{1,2,3\}$).

We now recall that the dual graph can be associated with the graph $\gamma\in\Gamma(g)$. All its faces will be triangles, except for one selected one, which is an $n$-gon. Therefore, we can rewrite the previous expansion as a sum over triangulations of a surface of genus $g$ with the boundary which is the $n$-gon:

\begin{equation}\label{Wick_triang}
    \left\langle H_{i_1i_2}H_{i_2i_3}\dots H_{i_ni_1} (W(H))^k\right\rangle_0=\sum\limits_{g=0}^{\infty}\sum\limits_{\{\tau\in \tilde{\Gamma}_n(g)\}}\sum\limits_{i_v\in I^2_N,\:v\in V(\tau)}\prod\limits_{t\in T(\tau)}\omega_t^{i_{v_1(t)}i_{v_2(t)}i_{v_3(t)}}
\end{equation}

where $\tilde{\Gamma}_n(g)$ is the set of triangulations of a surface of genus $g$ with boundary (on which $n$ triangulation vertices lie), $V(\tau)$ is the set of $\tau$ vertices, $T(\tau)$ is the set of triangles in $\tau$, $v_{j}(t)$ are vertices of triangle $t$.

The decomposition of \eqref{Wick_triang} does not depend on the explicit form of $\omega_{ijk}$ and on the index set $I^2_N$. The fact that the indices belong to a bounded lattice gives us a mapping of the vertices of each triangulation $\tau$ with fixed indices $i_v$ to the lattice $I^2_N$. This mapping extends linearly to both edges and vertices. Therefore, the sum $\sum\limits_{i_v}$ has the meaning of the sum over all triangulation mappings in $I^2_N$, such that the vertices of the boundary map bijectively and order-preserving into the set $\{i_1,\dots,i_n\}\subset I^2_N$.

Now let us fix the particular formula for $\omega_{ijk}$ as in \eqref{omega}. An important property here is that only triangulations make a nonzero contribution, in which, when mapped onto a lattice, all triangles are oriented consistently and are non-degenerate ($\epsilon^{ab}(i_a-j_a)(k_b-j_b)>0$). What this really means is that the triangulation mapping has to be a covering. This implies, firstly, that the image of a mapping can only be the interior of the region bounded by the polygon $i_1i_2\dots i_n$. Secondly, of the entire sum $\sum_{g}$ only the genus $g=0$ is non-zero, that is, the covering by the disk. Finally, the sum over the graphs is finite, since the number of vertices in a triangulation is limited due to the non-degeneracy of triangles (and, consequently, the surjectivity of the mapping of the set of vertices, since the covering is univalent and has no branches).

We note here that if we chose the vertices $i_1,\dots,i_n$ so that they do not form the boundary of the polygon, but, for example, are a polyline that wraps around zero twice, then branched coverings with non-surjective mapping of vertices could also be admissible.

Thus, we obtain the sum over all triangulations of the $\tau$ of $n$-gon $i_1\dots i_n$ with weights:

\begin{equation}
    \prod\limits_{t\in T(\tau)}e^{2\textrm{Vol}(t)(x^{v_1(t)}+x^{v_2(t)}+x^{v_3(t)})},
\end{equation}

where $2\textrm{Vol}(t)$ is twice the volume of the triangle $t$ arising from the factor $\epsilon^{ab}(i_a-j_a)(k_b-j_b)$. The product over all triangles can now be raised to the exponent and get the sum, which, after reducing similar terms, has the form:

\begin{equation}
    2\sum\limits_{a\in V(\tau)}x^a\sum\limits_{t\in T(\tau):\: a\in t}\textrm{Vol}(t)=2\sum\limits_{a\in V(\tau)}x^a\phi_a(\tau).
\end{equation}

Gathering everything together now, we obtain the required statement.

To prove the corollary, we note that the set of vertices in any maximal triangulation coincides with the set of interior points of the polygon and its vertices. Further, knowing the Euler characteristic of a surface with a boundary and the relation between the number of edges and triangles in a triangulation (3/2 of the number of triangles is equal to the number of edges + a correction to this formula taking into account the boundary), one can easily relate the number of interior points of the polygon $\ell^*(\Delta)$ with the number of triangles $|\tau|$.

\section{Appendix B. Virasoro-like algebra of trees}

Let us describe in more detail the structure of the potential of the tensor model \eqref{Tensor_model} for polytopes of dimension $d$. It has the following form

\begin{equation}
    \sum\limits_{n=1}^{\infty}\sum\limits_{\alpha\in \mathcal{T}_n}t_{(n,\alpha)} W^{n}_{\alpha}(T)
\end{equation}

\noindent where we grade the contributions to the potential with respect to the powers $(d-1)n$ of the tensor $T$. We define sets $\mathcal{T}_n$ and, accordingly, the potentials $W^n_{\alpha}(T)$ by recursion relations.

Let's start with the definition of $\mathcal{T}_1$. This is a set of one element corresponding to a $d$-dimensional simplex. The corresponding potential has already been written above \eqref{qubic_term}:

\begin{equation}
    W^1(T)=\sum\limits_{i_0\dots i_{d}}\Omega^{d+1}_{i_0\dots i_d}\prod\limits_{k=0}^{d}T_{i_0\dots\hat{i}_k\dots i_d}
\end{equation}
where $\Omega^{d+1}_{i_0\dots i_d}$ defined as

\begin{equation}
    \Omega^{d+1}_{i_0\dots i_d}\sim\begin{cases}
    1, \text{ if } \epsilon^{m_1\dots m_d}((i_1)_{m_1}-(i_0)_{m_1})\dots((i_d)_{m_d}-(i_0)_{m_d})>0 \\
    0, \text{ if } \epsilon^{m_1\dots m_d}((i_1)_{m_1}-(i_0)_{m_1})\dots((i_d)_{m_d}-(i_0)_{m_d})\leq 0
    \end{cases}.
\end{equation}

Note that the index set of $\Omega^{d+1}$ has the following structure:
\begin{itemize}
    \item given $d$-element subsets $\{i_0,\dots,\hat{i}_k,\dots,i_n\}$ corresponding to $T_{i_0\dots\hat{i}_k\dots i_d}$ tensors in potential $W^1(T)$,
    \item any two such subsets intersect in a $d-1$-element subset.
\end{itemize}
It is obvious that this corresponds to the structure of the simplicial complex of the boundary of a $d$-dimensional simplex (if we add all possible subsets to the $d$-element subsets). In higher powers of the potential in $T$, we want to have a similar structure of the simplicial complex of the boundary of some polyhedron, which we will denote by $\Delta(\alpha)$. To each such $d$-element subset we associate one factor $T$ in the corresponding term of the potential.

Now let the set $\mathcal{T}_{n}$ and the potentials $W^n_{\alpha},\:\alpha\in\mathcal{T}_n$ with the corresponding structure on $\Omega^n_{\ alpha}$ are given. Then the set $\mathcal{T}_{n+1}$ consists of all $d$-element subsets of all $\Omega^n$. In turn, each $W^n$ of the form:

\begin{equation}
    W^n_{\alpha}(T)=\Omega_{i_1\dots i_{l_1}\dots i_{l_2}\dots i_{l_d}\dots i_{q}}(T_{\dots}\dots T_{i_{l_1}\dots i_{l_k}}\dots T_{\dots})
\end{equation}

together with the $d$-element subset $\{i_{l_1},\dots,i_{l_d}\}$ we will associate the potential

\begin{equation}
    W^{n+1}_{\alpha'}(T)=\Omega_{i_1\dots i_{l_1} i_{l_0}\dots i_{l_2}\dots i_{l_d}\dots i_{n}}'(T_{\dots}\dots \hat{T}_{i_{l_1}\dots i_{l_k}}\dots T_{\dots})\prod\limits_{j=1}^{d}T_{i_{l_0}\dots i_{l_j}\dots i_{l_d}}
\end{equation}
where $\Omega'$ is defined by the formula

\begin{equation}\label{Omega_recurrence}
    \Omega_{i_1\dots i_{l_1} i_{l_0}\dots i_{l_2}\dots i_{l_d}\dots i_{n}}'=\begin{cases}
    \Omega_{i_1\dots i_{l_1}\dots i_{l_2}\dots i_{l_d}\dots i_{n}}\Omega^{d+1}_{i_{l_0}\dots i_{l_d}},\text{ if }\Delta(i_1\dots i_{l_1} i_{l_0}\dots i_{l_2}\dots i_{l_d}\dots i_{n}) \text{ is convex}\\
    0, \text{otherwise}
    \end{cases}
\end{equation}

Essentially, the recurrence relation \eqref{Omega_recurrence} simply corresponds to gluing an additional $d$-simplex to the polytope corresponding to $W^n_{\alpha}(T)$.

The resulting potential can be described in other, perhaps more descriptive terms. We will associate a tree whose set of edges for each vertex is ordered with each $W^n_{\alpha}(T)$. Namely, we will consider trees with fixed valence $d+1$ and open edges (having a beginning but no the end), where outgoing edges at each vertex are numbered from 0 to $d$. We will do it in the following way. Let us associate the potential $W^1(T)$ with the simplest tree, which is the point from which $d+1$ open edges emanate. Now let the potential $W^n_{\alpha}$ be associated with a tree. The \eqref{Omega_recurrence} recursion relation is equivalent to adding the simplest tree and gluing one of its open edges to an already existing open edge. The number of open edges is equal to the degree of the corresponding term of the potential in $T$. Thus, such trees correspond one-to-one to the elements of the sets $\mathcal{T}_n$.

Therefore, the potential of the tensor model can be rewritten as the sum over the trees described above. Let's denote the set of all such trees as $\mathcal{T}$.

\textbf{Ward identities and Virasoro-like algebra}

Now, in order to obtain the Ward identities \eqref{tensor_Ward_identities}, we just need to insert an operator of the following form under the integral of our tensor model:

\begin{equation}\label{L_tensor_operator}
    \hat{L}_{(\alpha,v)}=W^n_{\tau}(T,\dots,\frac{\partial}{\partial T},\dots,T)
\end{equation}

where $\frac{\partial}{\partial T}$ is in the place corresponding to the open edge $v$. So operators \eqref{L_tensor_operator} labeled by pairs $(\alpha,v)$ of some tree and a marked point of this tree.

It is easy to show that these operators, when acting on the $\langle\Delta\rangle$ correlator, can be rewritten as:

\begin{equation}
    \hat{L}_{(\alpha,v)}=-\mu \frac{\partial}{\partial t_{\alpha}}+\sum\limits_{\alpha'\in \mathcal{T}}t_{\alpha'}\sum\limits_{\tilde{\alpha}\in (\alpha'\#(\alpha,v))}\frac{\partial}{\partial t_{\tilde{\alpha}}} + \hat{I}_{(\alpha,v)}
\end{equation}
where the operator $\hat{I}_{\alpha}$ acts on $\langle\Delta\rangle$ by changing $\Delta$ in the analogous way as it was for the matrix model. Geometrically it glue the boundary of $\Delta(\alpha)$ to the boundary of $\Delta$ along the $d-1$-simplex associated with $v$.

Commuting such operators, we obtain the relations \eqref{Virasoro_like}.


\newpage


\bibliographystyle{unsrt}
\bibliography{references}

\end{document}